\documentclass[]{amsart}
\usepackage{amsmath}
\usepackage{amsfonts}
\usepackage{amssymb}
\usepackage{amsbsy}
\usepackage{amscd}
\usepackage{amsthm}
\usepackage{enumitem}
\usepackage{color}
\usepackage{framed}
\usepackage{graphicx}

\definecolor{sidegreen}{rgb}{0,0.45,0}

\newtheorem{theorem}{Theorem}[section]
\newtheorem{lemma}[theorem]{Lemma}
\newtheorem{corollary}[theorem]{Corollary}

\newtheorem{definition}{Definition}[section]
\newtheorem{proposition}[theorem]{Proposition}

\theoremstyle{remark}
\newtheorem{remark}[theorem]{Remark}
\newtheorem*{remark*}{Remark}

\numberwithin{equation}{section}

\newcommand{\T}{\mathbb{T}}
\newcommand{\cT}{\mathcal T}

\newcommand{\R}{\mathbb{R}}

\newcommand{\N}{\mathbb{N}}

\newcommand{\cO}{\mathcal O}

\newcommand{\eps}{\varepsilon}

\renewcommand{\phi}{\varphi}

\begin{document}

\title[Higher-order asymptotic formulae for Anosov flows]{Higher-order asymptotic formulae for periodic orbits of 3-dimensional volume-preserving Anosov flows and the Foulon conjecture}

\author{Andrey Gogolev}
\thanks{The author was partially supported by NSF grant DMS-2247747.}

\address{Department of Mathematics, The Ohio State University, Columbus, OH
43210, USA}
\email{gogolyev.1@osu.edu}

\begin{abstract}
Let $X^t$ be a volume-preserving Anosov flow of class $C^7$ on a closed
$3$-manifold.  We consider a sequence of periodic orbits which shadow a
homoclinic loop of a fixed periodic orbit $\gamma$, and derive
higher-order asymptotic expansions for their periods $T_n$ and for the logarithms of their unstable multipliers $\log\lambda_n$.  We give geometric
interpretations of the coefficients in these expansions, in particular in terms of the Anosov cocycle (\`a la Hurder-Katok~\cite{HK})
and the longitudinal Anosov cocycle (\`a la Foulon-Hasselblatt~\cite{FH}), evaluated
on $\gamma$.  Following the idea of De Simoi-Leguil-Vinhage-Yang
~\cite{DSLVY}, we compare these asymptotics when volume is the measure of maximal
entropy. In this case, the equilibrium-state characterizations of volume and the measure of maximal entropy synchronize periods and exponents:
$\log\lambda_n=hT_n$. Matching the coefficients of the two asymptotic expansions forces the cohomology classes of both cocycles to vanish.
 Then the work of
Foulon-Hasselblatt on the longitudinal Anosov cocycle~\cite{FH}, together with Foulon's work on the contact case~\cite{Fou}, yields the Foulon Conjecture for flows of class $C^7$: if volume is the measure of maximal entropy,
then, up to a finite cover, $X^t$ is smoothly conjugate either to a
constant-roof suspension of a hyperbolic automorphism of the $2$-torus or
to a finite cover of the constant-speed geodesic flow on a closed hyperbolic
surface.
\end{abstract}

\maketitle

\section{Introduction}

\subsection*{Statement of the results}

Let $M$ be a closed $3$-dimensional manifold and let
$X^t\colon M\to M$ be an Anosov flow generated by a vector field $X$. That is,
there is a $DX^t$-invariant splitting
$ TM=E^s\oplus\R X\oplus E^u,$
where $E^s$ is uniformly contracted and $E^u$ is uniformly expanded. We will
write $E^{0s}=E^s\oplus\R X$ and $E^{0u}=E^u\oplus\R X$
for the weak stable and weak unstable distributions, respectively.

\begin{theorem}
\label{thm_main}
Let $X^t\colon M\to M$ be a $C^r$, $r\ge 7$, volume-preserving Anosov flow on
a closed connected $3$-dimensional manifold. Assume that the 
invariant volume $m$ is the measure of maximal entropy. Then there exists a
finite cover such that the lifted flow is
$C^{r_*}$ conjugate\footnote{Here $r_*=r$ if $r\notin\N$, and
$r_*=(r-1)+\operatorname{Lip}$ if $r\in\N$. We note that the proof still works if $X^t$ is $C^{6+\operatorname{Lip}}$, but the author doesn't know how to lower the regularity assumption any lower.} to one of the following
algebraic flows:
\begin{enumerate}
\item a constant-roof suspension of a hyperbolic automorphism of $\T^2$;
\item a finite cover of the constant-speed geodesic flow on a closed
hyperbolic surface.
\end{enumerate}
\end{theorem}

The first result in this direction is probably the entropy rigidity theorem of Katok: if the metric entropy of a geodesic flow on a surface of negative curvature equals to its topological entropy then the metric has constant negative curvature~\cite{Kat}.
Foulon handled the more general case of contact flows and posed the even more general case of 3-dimensional volume-preserving Anosov flows as a conjecture~\cite{Fou}. De Simoi-Leguil-Vinhage-Yang developed an approach to this conjecture via asymptotic formula~\cite{DSLVY}. However the expansion for periods in~\cite[Lemma 3.4]{DSLVY} missed a term which lead the authors to conclude that the Anosov cocycle is trivial based only on the first order formula. 
Here we carefully derive the higher order formula and prove the above theorem by further pursuing the beautiful approach of De Simoi-Leguil-Vinhage-Yang. The need for higher order terms leads to a stronger regularity assumption $C^7$ (compared to original $C^5$ assumption in~\cite{DSLVY}).

Specifying to the setting of suspension flows we obtain
\begin{corollary}
\label{cor_torus}
Let $f\colon\T^2\to\T^2$ be a $C^r$, $r\ge7$, area-preserving Anosov
diffeomorphism and let $\phi=-\log J^u f$.
If $\phi$ is cohomologous to a $C^7$ function, then $f$ is
$C^{r_*}$ conjugate to a hyperbolic linear automorphism of $\T^2$. 
\end{corollary}

By ``cohomologous to a $C^7$ function" we precisely mean that
$
\phi=\bar\phi+u\circ f-u,
$
where $\bar\phi$ is $C^7$ and $u$ is continuous.
\subsection*{Organization}

In Section~2 we explain some needed preliminaries. We recall how the entropy assumption synchronizes the periods and the Lyapunov exponents of periodic points. We introduce an elaborate version of Hurder-Katok adapted charts for Anosov flows, the Anosov cocycle and the longitudinal Anosov cocycle (or the Foulon-Hasselblatt cocycle) and recall their rigidity properties. Finally we also introduce a finite regularity version of the Birkhoff normal form at a periodic point. 

In Section~3 we prove the higher-order asymptotics for the periods of shadowing orbits. And in Section~4 we do the same for their Lyapunov exponents.  Finally, in Section~5 we deduce our main result and the corollary.

\subsection*{Comments} The author is very grateful to Martin Leguil from whom he learnt about asymptotic formulae for the periods in July 2022. At the same time the author concluded that the main decreasing term of the order $n\mu^n$ was missing in the formula for the periods and communicated this error to M. Leguil. It was clear that pursuing higher order asymptotics is a plausible route for obtaining Foulon conjecture. The author began actively thinking about it in November 2025 after learning about~\cite{EVY} since it provided such a clean normal form at a periodic orbit.\footnote{While psychologically the normal form of~\cite{EVY} was very important, the proof here only relies on the Birkhoff normal form.} 
 By May 2026 the author was convinced that derivation of higher order asymptotics is possible and it was the matter of getting precise value of the coefficients in the formula which appeared to be a very tedious calculus task. At this point AI help (first Claude Opus 4.8 and then ChatGPT Sol 5.6) became invaluable for doing calculations. The most challenging one was the polynomial drift expansion (Lemma~\ref{lemma_cross_form}) which the author could derive, but not  with the needed $C^1$-control. In August 2026 ChatGPT pointed out a paper of M. Gonchenko~\cite{Gon13} which had a lemma very close to what was needed and the task was finished soon after. The author also used Claude and ChatGPT to generate initial drafts on most of the proofs in this paper, all of which were substantially, but not completely, rewritten by the author. The author would like to thank OpenAI for providing Researcher Workspace.

The author would like to thank Giovanni Forni for feedback on early drafts.

\section{Preliminaries}
\label{sec_preliminaries}

We pass to the connected finite cover which orients $E^s$ and $E^u$.
The lifted volume is still the measure of maximal entropy, hence, it is enough to work with the
lifted flow. Also note that for the lifted flow all periodic orbits have positive multipliers.

\subsection{Synchronization}

\leavevmode
We begin with a well-known consequence of the assumption that
volume is the measure of maximal entropy. Fix a background Riemannian metric
on $M$. Let
$$
\phi^u(x)=\left.\frac{d}{dt}\right|_{t=0}\log J^uX^t(x)
$$
be the infinitesimal logarithmic unstable Jacobian. The volume $m$ is the SRB measure and,
hence, is the equilibrium state of the potential $-\phi^u$. The pressure of
this potential is zero~\cite[Theorem~3.3 and Proposition~4.4]{BR}. By assumption, the same measure is the equilibrium state of
the zero potential, whose pressure is the topological entropy
$h=h_{\mathrm{top}}(X^1)$.
It follows that the potentials are cohomologous, that is,
$$
\phi^u=h+Xu,
$$
where $u$ is a H\"older function differentiable in the flow
direction; see, e.g.,~\cite[Lemma~2.4]{Sh26}.
Integrating over a periodic orbit gives the following lemma.

\begin{lemma}
\label{lemma_sync}
If $p=X^T(p)$, then
\begin{equation*}
\label{eq_sync}
\log J^uX^T(p)=hT.
\end{equation*}
\end{lemma}

\subsection{Non-stationary Birkhoff normal form and Anosov cocycle}
\label{subsec_adapted_charts}

Let $X^t$ be a $C^r$, $r\ge7$,
volume-preserving Anosov flow on a closed $3$-dimensional manifold, and let
$\Omega$ be the invariant volume form.
Here we elaborate on the construction of adapted transverse coordinates by
Hurder and Katok~\cite[Sections~4]{HK}. We write the transition
maps through cubic order in non-stationary Birkhoff normal form, so that the
Anosov cocycle can be read directly from this normal form.

\begin{proposition}
\label{prop_nonstationary_birkhoff}
There exist $\eps>0$ and a family of embeddings
$$
\Psi_q\colon(-\eps,\eps)^2\longrightarrow M,
\qquad \Psi_q(0,0)=q,  \,\,\,\,q\in M,
$$
with the following properties.
\begin{enumerate}
\item The embeddings are $C^{r-1}$ and
$q\mapsto\Psi_q$ is H\"older continuous in the $C^{r-1}$ topology. Their
images $\cT_q$ are transverse to the flow.
\item The curves $\Psi_q(\cdot,0)$ and $\Psi_q(0,\cdot)$ parametrize the
local strong stable and local strong unstable manifolds of $q$, respectively.
\item Adapted coordinates respect area froms
$$
\Psi_q^*(\iota_X\Omega)=dx\wedge dy.
$$
\item Let $F_{q,t}\colon\cT_q\to\cT_{X^t(q)}$ be the local Poincar\'e map and consider its coordinate representation
$$
\mathcal F_{q,t}=\Psi_{X^t(q)}^{-1}\circ F_{q,t}\circ\Psi_q.
$$
Then,
\begin{multline}
\label{eq_adapted_return}
\mathcal F_{q,t}(x,y)=
\big(\mu(q,t)x\big(1-A(q,t)xy\big),
\\
\mu(q,t)^{-1}y\big(1+A(q,t)xy\big)\big)+R_{q,t}(x,y),
\end{multline}
where $\mu(q,t)>0$. Moreover, $R_{q,t}$ vanishes to order three at the
origin:
$$
D^jR_{q,t}(0,0)=0,\qquad 0\leq j\leq 3.
$$
\end{enumerate}
\end{proposition}


\begin{definition}
The coefficient $A(q,t)$ in~\eqref{eq_adapted_return} is the
\emph{Anosov cocycle}. Equivalently,
$$
A(q,t)=-\frac{1}{2\mu(q,t)}
\partial_x^2\partial_y(\mathcal F_{q,t})_1(0,0)
=\frac{\mu(q,t)}2
\partial_x\partial_y^2(\mathcal F_{q,t})_2(0,0).
$$
\end{definition}

It is routine to check that $A$ is a H\"older continuous additive cocycle.
Hurder and Katok~\cite[Theorem~3.4]{HK} classified volume-preserving Anosov
flows with smooth weak unstable distribution as precisely those for which
$[A]=0$. 
However, we will not use this result in this paper.

\begin{remark}
If $p=X^{T_0}(p)$, then $\mathcal F_{p,T_0}$ is an area-preserving germ at a
hyperbolic fixed point. Recall the standard Birkhoff normal form
$$
(x,y)\longmapsto
\big(\mu xB(xy),\mu^{-1}yB(xy)^{-1}\big),
\qquad B(0)=1,
$$
formula~\eqref{eq_adapted_return} gives
$$
B'(0)=-A(p,T_0).
$$
Thus $A(p,T_0)$ is the hyperbolic twist. See Birkhoff~\cite[Chapter~VIII, Section~1,
pp.~211--213]{Bir27} and Hurder-Katok~\cite[Section~5]{HK}.
\end{remark}

%


\subsection{Longitudinal Anosov cocycle}
\label{subsec_cocycles}

For $q\in M$ and $T>0$, consider the adapted transversals based at the
 points $q$ and $X^T(q)$. We write the transition time between them as
$T+\tau_{q,T}$; thus $\tau_{q,T}$ is the unique $C^{r-1}$ function
such that
$$
X^{T+\tau_{q,T}(x,y)}\Psi_q(x,y)
=\Psi_{X^T(q)}\big(\mathcal F_{q,T}(x,y)\big),
\qquad \tau_{q,T}(0,0)=0.
$$
Since the coordinate axes parametrize the strong stable and strong unstable
manifolds, we have
$$
\tau_{q,T}(x,0)=\tau_{q,T}(0,y)=0.
$$

\begin{definition}
The \emph{longitudinal Anosov cocycle}, or the
\emph{Foulon-Hasselblatt cocycle}, is defined as
\begin{equation*}
\label{eq_longitudinal_cocycle}
K(q,T)=\partial_x\partial_y\tau_{q,T}(0,0).
\end{equation*}
\end{definition}

Then $K$ is a cocycle whose cohomology class $[K]$ is independent of the choice of adapted charts~\cite{FH}.

\begin{theorem}[{Foulon-Hasselblatt~\cite[Theorem~3]{FH}}]
\label{thm_Foulon_Hasselblatt}
Let $X^t$ be a $C^r$, $r\ge2$, volume-preserving Anosov flow on a closed
$3$-manifold. Then
$$
\begin{aligned}
{}[K]=0
&\quad\Longleftrightarrow\quad E^s\oplus E^u\in C^{r-1}\\
&\quad\Longleftrightarrow\quad
X^t\text{ is a contact flow or a constant-roof suspension}.
\end{aligned}
$$
\end{theorem}


\subsection{A finite-order Birkhoff normal form}
\label{subsec_finite_birkhoff}

Fix a periodic orbit $\gamma$, a point $p\in\gamma$, and let $T_0$ be its
period. The adapted chart of
Proposition~\ref{prop_nonstationary_birkhoff} gives $C^{r-1}$
area-preserving coordinates $(\xi,\eta)$ on $\cT_p$ relative to which the
return map has the form
\begin{equation}
\label{eq_first_birkhoff}
F(\xi,\eta)=\big(\mu\xi(1-a\xi\eta)+R_1(\xi,\eta),
\mu^{-1}\eta(1+a\xi\eta)+R_2(\xi,\eta)\big),
\end{equation}
where $\mu=\mu(p,T_0)\in(0,1)$, $a=A(p,T_0)$, and
$D^jR_1(0,0)=D^jR_2(0,0)=0$ for $0\leq j\leq 3$.
We will need more
refined coordinates at the periodic point $p$.

The following proposition is essentially the same as a more general higher order normal form~\cite[Lemma~1.1]{Gon13}. However, we save one extra derivative (compared to $C^{r-4}$) and we have a better control on the error terms. We include a proof in Appendix~\ref{app_normal_forms}.

\begin{proposition}
\label{prop_normal_form}
There exist $C^{r-3}$ area-preserving coordinates $(x,y)$ on $\cT_p$ in which
\begin{multline}
\label{eq_finite_moser}
F(x,y)=\big(\mu x\big(1-axy+b(xy)^2\big)+r_1(x,y),
\\
\mu^{-1}y\big(1+axy+(a^2-b)(xy)^2\big)+r_2(x,y)\big),
\end{multline}
where $0<\mu<1$ is the stable multiplier of $p$,
$a=A(p,T_0)$, and the error terms have the form
\begin{equation}
\label{eq_moser_remainder}
\begin{split}
r_1(x,y)&=x^3y^2\rho_1(x,y),\\
r_2(x,y)&=x^2y^3\rho_2(x,y),
\end{split}
\end{equation}
where $\rho_1,\rho_2\in C^{r-6}$ and $\rho_1(0,0)=\rho_2(0,0)=0$.

Let $T_0+\tau$ be the return time in these coordinates. Then there exists
$v\in C^{r-5}$, vanishing on both coordinate axes, such that
$\widehat\tau=\tau+v\circ F-v$ has the expansion
\begin{equation}
\label{eq_finite_roof}
\widehat\tau(x,y)=K(p,T_0)xy+O((xy)^2).
\end{equation}
\end{proposition}


\section{Asymptotic formula for the periods}
\label{sec_period_formula}

The formula below refines the first-order expansion derived
in~\cite{GRH}. We fix a periodic orbit $\gamma$ and the Birkhoff coordinates
given by Proposition~\ref{prop_normal_form}.

Let $\cO$ be an orbit homoclinic to $\gamma$. Let
$p_{\rm out}\in\cO\cap W^u_{\rm loc}(p)\cap\cT_p$ be its last intersection with the transversal $\cT_p$ 
before the excursion and let
$p_{\rm in}\in\cO\cap W^s_{\rm loc}(p)\cap\cT_p$ be its first intersection
after the excursion. After replacing these points by
$F^{-j}(p_{\rm out})$ and $F^j(p_{\rm in})$ for some $j\ge 0$, and making $\cT_p$ smaller if needed, we may
assume that the open excursion from $p_{\rm out}=(0, y_{\rm out})$ to $p_{\rm in}=(x_{\rm in}, 0)$ is disjoint
from $\cT_p$. The excursion defines a Poincar\'e map $G$ from a neighborhood
of $p_{\rm out}$ in $\cT_p$ to a neighborhood of $p_{\rm in}$ in $\cT_p$ with $G(p_{\rm out})=p_{\rm in}$, and
an excursion-time function $T_{\rm exc}$ on its domain, which is characterized by
$$
X^{T_{\rm exc}(q)}(q)=G(q).
$$
In the refined Birkhoff coordinates, $G$ and $T_{\rm exc}$ are $C^{r-3}$.

For all sufficiently large $n$, the associated pseudo-orbit is shadowed by a
periodic orbit $\gamma_n$ which intersects $\cT_p$ $(n+1)$ times. Label its intersections
$p_i=(x_i,y_i)\in\cT_p$, $0\leq i\leq n$, ordered by time so that
$p_0$ is close to $p_{\rm in}$ and $p_n$ is close to $p_{\rm out}$. Then $p_n$ belongs to
the domain of $G$ for all sufficiently large $n$. We have
$$
p_i=F^i(p_0),\qquad p_0=G(p_n).
$$
We denote the period of $\gamma_n$ by $T_n$.

Throughout all local calculations that follow we work in the first quadrant:
$x_i>0$ and $y_i>0$ for $0\leq i\leq n$. Indeed, this is just a matter of properly orienting the local coordinates. (Recall that $p$ has positive multipliers after passing to a finite cover.)

\begin{proposition}
\label{prop_period_formula}
There exist constants $C$ and $c\neq0$ such that
\begin{equation*}
\label{eq_higher_period}
\begin{split}
T_n={}&nT_0+T_{\rm exc}(p_{\rm out})+C\mu^n
+cK(p,T_0)n\mu^n\\
&-c^2A(p,T_0)K(p,T_0)n^2\mu^{2n}+O(n\mu^{2n}).
\end{split}
\end{equation*}
\end{proposition}

For the sake of doing asymptotics we subdivide $\gamma_n$ into the {\it excursion part} from $p_n$ to $p_0$ which is governed by $G$ and return time $T_{\rm exc}$ and the {\it local passage} from $p_0$ to $p_n$ which is governed by the local dynamics $F$ and the roof function $T_0+\tau$. 
\begin{equation}
\label{eq_split}
T_n=\sum_{i=1}^n(T_0+\tau(x_i,y_i))+T_{\rm exc}(p_n)
\end{equation}

We will first expand the product $I_n=x_0y_0$ in terms of $\mu^n$ and then
use it to expand $T_n$.

It will be useful to consider the local passage $(p_i)_{i=0}^n$ as parametrized via the $x$-coordinate of the first point $x_0\approx x_{\rm in}$ and the $y$-coordinate of the last point $y_n\approx y_{\rm out}$. In particular, we will view $x_n$ and $y_0$ as functions of $(x_0, y_n)$.
The following lemma concerns only the local dynamics $F$ and gives the
$(x_0,y_n)$-parametrization and asymptotic expansion of $x_n$ and $y_0$.\footnote{Note that while not explicitly indicated in the notation $y_0$ also depends on $n$ as well as the whole local passage $(p_i)_{i=0}^n$.}
Then Lemma~\ref{lemma_In_estimate} incorporates the
excursion map $G$ and gives asymptotics of $I_n$.

\begin{remark*}
The next lemma is almost the same as a lemma by M. Gonchenko~\cite[Lemma~1.2]{Gon13}. We
follow contracting summation-operator proof from
\cite[Section~4.2.2]{Gon13}. It is stated in~\cite[Lemma 1.2]{Gon13} that $C^1$ control can be achieved, but the details were not explicitly given. We carefully calculate the coefficients the polynomial $\Pi_n$ and establish $C^1$ control of the remainder. 
\end{remark*}

\begin{lemma}
\label{lemma_cross_form}
For all sufficiently large $n$, we have
\begin{equation}
\label{eq_cross_form}
\begin{split}
x_n&=\mu^nx_0\Pi_n(s)+\Delta_{1,n}(x_0,y_n),\\
y_0&=\mu^ny_n\Pi_n(s)+\Delta_{2,n}(x_0,y_n),
\end{split}
\end{equation}
where $s=\mu^nx_0y_n$ and the main term is given by the polynomial
\begin{equation}
\label{eq_passage_polynomial}
\Pi_n(s)=1-nas+
\left(\frac32a^2n^2+\left(b-\frac12a^2\right)n\right)s^2,
\end{equation}
and the errors are given by functions $\Delta_{1,n}$ and $\Delta_{2,n}$ which are defined on a neighborhood of $(x_{\rm in}, y_{\rm out})$ and are controlled in $C^1$ norm
$$
\|\Delta_{1,n}\|_{C^1}+\|\Delta_{2,n}\|_{C^1}\leq C\mu^{3n}.
$$
\end{lemma}

\begin{remark*} The last summand of the main term $\mu^nx_0\Pi_n(s)$ is of order $\mu^{3n}$ in the $C^0$-norm, but of order $\mu^{2n}$ in the $C^1$-norm.
\end{remark*}

\begin{proof}
$\quad$\\
\noindent{\bfseries Step 1: }{\it Local passage as a fixed point of a
summation operator.}

Let
$$
B(u)=1-au+bu^2,
$$
where $a$ and $b$ are the constants coming from the Birkhoff normal form~\eqref{eq_finite_moser}.
Define the integrable map
$$
F_0(x,y)=\big(\mu xB(xy),\mu^{-1}yB(xy)^{-1}\big).
$$
Then~\eqref{eq_finite_moser} can be written as
\begin{equation}
\label{eq_integrable_remainder}
F=F_0+\mathcal R,
\qquad
\mathcal R(x,y)=\big(x^3y^2\rho_1(x,y),x^2y^3\rho_2(x,y)\big),
\end{equation}
where the two functions $\rho_j$ are $C^1$ and vanish at the origin. Here we
used
$$
B(u)^{-1}=1+au+(a^2-b)u^2+O(u^3)
$$
and absorbed the $O(u^3)$ term into the second remainder.

We rescale the coordinates $(x_i,y_i)$, $0\leq i\leq n$, to be of order one:
$$
\widehat x_i=\mu^{-i}x_i,\qquad
\widehat y_i=\mu^{i-n}y_i.
$$
Then 
$x_iy_i=\mu^n\widehat x_i\widehat y_i$.
The decomposition~\eqref{eq_integrable_remainder} then reads
\begin{align}
\widehat x_{i+1}
&=\widehat x_iB(\mu^n\widehat x_i\widehat y_i)
+\mathcal E_{i,n}(\widehat x_i,\widehat y_i),
\label{eq_scaled_x}\\
\widehat y_i
&=\widehat y_{i+1}B(\mu^n\widehat x_i\widehat y_i)
+\mathcal H_{i,n}(\widehat x_i,\widehat y_i),
\label{eq_scaled_y}
\end{align}
where
\begin{align*}
\mathcal E_{i,n}(x,y)
&=\mu^{2n-1}x^3y^2\rho_1(\mu^ix,\mu^{n-i}y),\\
\mathcal H_{i,n}(x,y)
&=-\mu^{2n+1}B(\mu^nxy)x^2y^3
\rho_2(\mu^ix,\mu^{n-i}y).
\end{align*}
Recall that $\rho_1$ and $\rho_2$ are at least $C^1$. Vanishing of $\rho_1$ and $\rho_2$ at the origin implies
\begin{equation*}
\label{eq_scaled_error}
\|\mathcal E_{i,n}\|_{C^1}+\|\mathcal H_{i,n}\|_{C^1}
\leq C\mu^{2n}(\mu^i+\mu^{n-i}).
\end{equation*}
Summing gives
\begin{equation}
\label{eq_summed_scaled_error}
\sum_{i=0}^{n-1}
\big(\|\mathcal E_{i,n}\|_{C^1}
+\|\mathcal H_{i,n}\|_{C^1}\big)\leq C\mu^{2n}.
\end{equation}

We view~\eqref{eq_scaled_x} and~\eqref{eq_scaled_y} as a system of equations
with $\widehat x_0=x_0$ and $\widehat y_n=y_n$ given. We now reinterpret this
system as a fixed point problem for a certain summation operator. For a
sequence $(\xi,\eta)=((\xi_i,\eta_i))_{i=0}^n$, let
$$
\|(\xi,\eta)\|_\infty
=\max_{0\leq i\leq n}(|\xi_i|+|\eta_i|),\qquad
\mathcal B_n(R)=\{(\xi,\eta):\|(\xi,\eta)\|_\infty\leq R\}.
$$
We equip $\mathcal B_n(R)$ with the uniform metric
$$
d_\infty((\xi,\eta),(\widetilde\xi,\widetilde\eta))
=\|(\xi,\eta)-(\widetilde\xi,\widetilde\eta)\|_\infty.
$$
For each $n$, $\mathcal B_n(R)$ is a closed ball in $(\R^2)^{n+1}$.
Choose $R$, independently of $n$, so that
$|x_0|+|y_n|<R/2$. Writing
\eqref{eq_scaled_x} and~\eqref{eq_scaled_y} in difference form we have
\begin{align*}
\widehat x_{i+1}-\widehat x_i
&=\widehat x_i(B(\mu^n\widehat x_i\widehat y_i)-1)
+\mathcal E_{i,n}(\widehat x_i,\widehat y_i),\\
\widehat y_i-\widehat y_{i+1}
&=\widehat y_{i+1}(B(\mu^n\widehat x_i\widehat y_i)-1)
+\mathcal H_{i,n}(\widehat x_i,\widehat y_i).
\end{align*}
Summing the first equality from $0$ to $i-1$ and the second from $i$ to
$n-1$ defines the summation operator
$\mathcal S_n(\xi,\eta)=(\xi',\eta')$ by
\begin{align*}
\xi_i'={}&\widehat x_0+\sum_{j=0}^{i-1}
\big(\xi_j(B(\mu^n\xi_j\eta_j)-1)
+\mathcal E_{j,n}(\xi_j,\eta_j)\big),\\
\eta_i'={}&\widehat y_n+\sum_{j=i}^{n-1}
\big(\eta_{j+1}(B(\mu^n\xi_j\eta_j)-1)
+\mathcal H_{j,n}(\xi_j,\eta_j)\big),
\end{align*}
whose fixed points are precisely the solutions
of~\eqref{eq_scaled_x}--\eqref{eq_scaled_y} with
$\widehat x_0=x_0$ and $\widehat y_n=y_n$.

On $\mathcal B_n(R)$, the $C^1$ norms of the functions
$$
(x,y)\mapsto x(B(\mu^nxy)-1)\quad\text{and}\quad
(x,y,z)\mapsto z(B(\mu^nxy)-1)
$$
are bounded by $C\mu^n$. Together with~\eqref{eq_summed_scaled_error}, this
gives
$$
\|\mathcal S_n(\xi,\eta)-\mathcal S_n(\widetilde\xi,
\widetilde\eta)\|_\infty
\leq (Cn\mu^n+C\mu^{2n})
\|(\xi,\eta)-(\widetilde\xi,\widetilde\eta)\|_\infty.
$$
The same estimates give
$$
\|\mathcal S_n(\xi,\eta)\|_\infty
\leq x_0+y_n+Cn\mu^n+C\mu^{2n}<R
$$
for all sufficiently large $n$. Thus $\mathcal S_n$ preserves
$\mathcal B_n(R)$ and is a contraction towards the fixed point with Lipschitz constant $1/2$ with respect to $\|\cdot\|_\infty$ norm.

\smallskip
\noindent{\bfseries Step 2: }{\it Truncated summation operator and estimates
on the distance between fixed points.}

The rescaled dynamics of $F_0$ is
\begin{align}
\widehat x_{i+1}^0
&=\widehat x_i^0B(\mu^n\widehat x_i^0\widehat y_i^0),
\label{eq_scaled_x_zero}\\
\widehat y_i^0
&=\widehat y_{i+1}^0B(\mu^n\widehat x_i^0\widehat y_i^0),
\label{eq_scaled_y_zero}
\end{align}
for $0\leq i\leq n-1$, with 
$\widehat x_0^0=x_0$ and $\widehat y_n^0=y_n$. Note that $\widehat x_{i+1}^0\widehat y_{i+1}^0=\widehat x_i^0\widehat y_i^0$.
Thus the product $\widehat x_i^0\widehat y_i^0$ is independent of $i$.
As before, writing~\eqref{eq_scaled_x_zero} and~\eqref{eq_scaled_y_zero} in difference
form and summing up defines the $F_0$ summation operator
$\mathcal S_n^0$. Explicitly, if
$\mathcal S_n^0(\xi,\eta)=(\xi',\eta')$, then
\begin{align*}
\xi_i'={}&\widehat x_0+\sum_{j=0}^{i-1}
\xi_j(B(\mu^n\xi_j\eta_j)-1),\\
\eta_i'={}&\widehat y_n+\sum_{j=i}^{n-1}
\eta_{j+1}(B(\mu^n\xi_j\eta_j)-1).
\end{align*}
Thus the fixed point of $\mathcal S_n^0$ is precisely the rescaled local
passage of $F_0$ satisfying the prescribed ``boundary condition."
The same argument shows that $\mathcal S_n^0$ preserves
$\mathcal B_n(R)$ and is a contraction with Lipschitz constant $1/2$ and
 $(\widehat x_i^0,\widehat y_i^0)$ is its fixed point.  By~\eqref{eq_summed_scaled_error},
\begin{equation}
\label{eq_operator_comparison}
\sup_{\mathcal B_n(R)}\|\mathcal S_n-\mathcal S_n^0\|_\infty
+\sup_{\mathcal B_n(R)}\|D\mathcal S_n-D\mathcal S_n^0\|_\infty
\leq C\mu^{2n}.
\end{equation}
Here and below, $D\mathcal S_n$ and $D\mathcal S_n^0$ denote derivatives in
the sequence variables. 

Recall that for each fixed $(x_0,y_n)$, the local passages
$(\widehat x,\widehat y)$ and $(\widehat x^0,\widehat y^0)$ are the fixed points of $\mathcal S_n$ and $\mathcal S_n^0$, respectively. Then
\begin{align*}
\|(\widehat x,\widehat y)-(\widehat x^0,\widehat y^0)\|_\infty
&\leq\|\mathcal S_n(\widehat x,\widehat y)
-\mathcal S_n(\widehat x^0,\widehat y^0)\|_\infty
+\|(\mathcal S_n-\mathcal S_n^0)(\widehat x^0,\widehat y^0)\|_\infty\\
&\leq \frac12\|(\widehat x,\widehat y)
-(\widehat x^0,\widehat y^0)\|_\infty
+\sum_{i=0}^{n-1}\big(\|\mathcal E_{i,n}\|_{C^0}
+\|\mathcal H_{i,n}\|_{C^0}\big).
\end{align*}
Thus~\eqref{eq_operator_comparison} gives
$$
\|(\widehat x,\widehat y)-(\widehat x^0,\widehat y^0)\|_\infty
=\max_{0\leq i\leq n}
\big(|\widehat x_i-\widehat x_i^0|
+|\widehat y_i-\widehat y_i^0|\big)
\leq 2C\mu^{2n}.
$$

Having estimated the uniform norm, we now estimate the $C^1$ norms as functions of $(x_0,y_n)$. The parameter-dependent contraction mapping theorem implies
that the fixed points are $C^1$ functions of $(x_0,y_n)$. Displaying the
parameter dependence, the fixed point equation for $\mathcal S_n$ is
$$
(\widehat x,\widehat y)
=\mathcal S_n(\widehat x,\widehat y;x_0,y_n).
$$
Differentiating using the chain rule gives
\begin{align*}
D_{(x_0,y_n)}(\widehat x,\widehat y)
={}&D\mathcal S_n(\widehat x,\widehat y)
D_{(x_0,y_n)}(\widehat x,\widehat y)
+\partial_{(x_0,y_n)}\mathcal S_n,
\end{align*}
and therefore
$$
\big(I-D\mathcal S_n(\widehat x,\widehat y)\big)
D_{(x_0,y_n)}(\widehat x,\widehat y)
=\partial_{(x_0,y_n)}\mathcal S_n.
$$
Here $D\mathcal S_n$ denotes the derivative in the sequence variables, while
$\partial_{(x_0,y_n)}\mathcal S_n$ is the differential with respect to 
$(x_0,y_n)$. By definition of $\mathcal S_n$, it is simply the map
$$
(h_1,h_2)\longmapsto
\big((h_1,h_1,\ldots,h_1),(h_2,h_2,\ldots,h_2)\big),
$$
which has norm one.

The same formulae
apply to $\mathcal S_n^0$. Since both contractions have Lipschitz constant 
 $1/2$, the operators on the left are invertible and their inverses have norm
at most $2$. The derivative of $\mathcal S_n^0$ is $Cn\mu^n$-Lipschitz on
$\mathcal B_n(R)$. Therefore, by the preceding $C^0$ estimate
and~\eqref{eq_operator_comparison},
\begin{align*}
&\|D\mathcal S_n(\widehat x,\widehat y)
-D\mathcal S_n^0(\widehat x^0,\widehat y^0)\|\\
&\qquad\leq
\|D\mathcal S_n(\widehat x,\widehat y)
-D\mathcal S_n^0(\widehat x,\widehat y)\|
+Cn\mu^n\|(\widehat x,\widehat y)
-(\widehat x^0,\widehat y^0)\|_\infty
\leq C\mu^{2n}.
\end{align*}
Note that
$\partial_{(x_0,y_n)}\mathcal S_n
=\partial_{(x_0,y_n)}\mathcal S_n^0$.
Subtracting the differentiated fixed point equations gives
\begin{align*}
0
&=\big(I-D\mathcal S_n\big)
D_{(x_0,y_n)}(\widehat x,\widehat y)-
\big(I-D\mathcal S_n^0\big)
D_{(x_0,y_n)}(\widehat x^0,\widehat y^0)\\
&=\big(I-D\mathcal S_n\big)
D_{(x_0,y_n)}\big((\widehat x,\widehat y)
-(\widehat x^0,\widehat y^0)\big)-
\big(D\mathcal S_n-D\mathcal S_n^0\big)
D_{(x_0,y_n)}(\widehat x^0,\widehat y^0)\\
&=\big(I-D\mathcal S_n\big)
\Big[D_{(x_0,y_n)}\big((\widehat x,\widehat y)
-(\widehat x^0,\widehat y^0)\big)\\
&\qquad-
\big(I-D\mathcal S_n\big)^{-1}
\big(D\mathcal S_n-D\mathcal S_n^0\big)
D_{(x_0,y_n)}(\widehat x^0,\widehat y^0)\Big],
\end{align*}
where $D\mathcal S_n$ and $D\mathcal S_n^0$ are evaluated
at $(\widehat x,\widehat y)$ and $(\widehat x^0,\widehat y^0)$, respectively.
Since $I-D\mathcal S_n$ is invertible, the expression
in brackets vanishes.
The second differentiated fixed point equation and the inverse bound give
$$
\|D_{(x_0,y_n)}(\widehat x^0,\widehat y^0)\|\leq2.
$$
Consequently,
\begin{align*}
&\big\|D_{(x_0,y_n)}\big((\widehat x,\widehat y)
-(\widehat x^0,\widehat y^0)\big)\big\|_\infty\\
&\leq\big\|\big(I-D\mathcal S_n\big)^{-1}\big\|_\infty
\big\|D\mathcal S_n-D\mathcal S_n^0\big\|_\infty
\big\|D_{(x_0,y_n)}(\widehat x^0,\widehat y^0)\big\|_\infty
\leq4C\mu^{2n}.
\end{align*}
Hence,
\begin{equation*}
\label{eq_scaled_comparison}
\max_{0\leq i\leq n}
\big(\|\widehat x_i-\widehat x_i^0\|_{C^1}
+\|\widehat y_i-\widehat y_i^0\|_{C^1}\big)\leq4C\mu^{2n}.
\end{equation*}
Here, again, we view the local passages $(\widehat x_i,\widehat y_i)$ and
$(\widehat x_i^0,\widehat y_i^0)$ as functions of $(x_0,y_n)$.
Since $x_n=\mu^n\widehat x_n$ and $y_0=\mu^n\widehat y_0$, this implies
\begin{equation}
\label{eq_cross_F_F0_comparison}
\|x_n-x_n^0\|_{C^1}+\|y_0-y_0^0\|_{C^1}
\leq C\mu^{3n}.
\end{equation}

\smallskip
\noindent{\bfseries Step 3: }{\it Asymptotics of the truncated fixed point.}

In this step we derive asymptotics for the $F_0$-passage.
Set $u=\mu^n\widehat x_i^0\widehat y_i^0$, which is independent of $i$.
Since $\widehat x_0^0=x_0$ and $\widehat y_n^0=y_n$, we have
$$
\widehat x_n^0=x_0B(u)^n,\qquad
\widehat y_0^0=y_nB(u)^n,\qquad u=sB(u)^n,
$$
where $s=\mu^nx_0y_n$. Note that $u=O(s)$. Since
$$
B(u)^n=1-nau+O(n^2u^2),
$$
the equation $u=sB(u)^n$ gives
$$
u=s-nasu+O(n^2su^2)=s+O(ns^2).
$$
Substituting this estimate back gives
\begin{align*}
u&=s-nas\big(s+O(ns^2)\big)+O(n^2s^3)\\
&=s-nas^2+O(n^2s^3).
\end{align*}
It follows that
\begin{align*}
\widehat x_n^0
&=x_0B(u)^n\\
&=x_0\left(1-nas+
\left(nb+\frac{n(3n-1)}2a^2\right)s^2
+O(n^3s^3)\right)\\
&=x_0\Pi_n(s)+O(n^3s^3).
\end{align*}
Analogously,
$$
\widehat y_0^0=y_n\Pi_n(s)+O(n^3s^3).
$$
The remainders and their first derivatives with respect to $s$ are bounded by
$Cn^3s^3$ and $Cn^3s^2$, respectively. Hence, after substituting
$s=\mu^nx_0y_n$, the identities
$\partial_{x_0}s=\mu^ny_n$ and $\partial_{y_n}s=\mu^nx_0$ show that the
derivatives of remainders with respect to $(x_0,y_n)$ are bounded by $Cn^3\mu^{3n}$.
Therefore,
$$
\max\left\{\|\widehat x_n^0-\widehat x_0\Pi_n(s)\|_{C^1},
\|\widehat y_0^0-\widehat y_n\Pi_n(s)\|_{C^1}\right\}
\leq Cn^3\mu^{3n}.
$$
Equivalently,
\begin{equation}
\label{eq_F0_polynomial_comparison}
\max\left\{\|x_n^0-\mu^nx_0\Pi_n(s)\|_{C^1},
\|y_0^0-\mu^ny_n\Pi_n(s)\|_{C^1}\right\}
\leq Cn^3\mu^{4n}.
\end{equation}

\smallskip
\noindent{\bfseries Step 4: }{\it Assembling everything together.}

The $C^1$ bounds~\eqref{eq_cross_F_F0_comparison}
and~\eqref{eq_F0_polynomial_comparison} give
\begin{align*}
&\max\left\{\|x_n-\mu^nx_0\Pi_n(s)\|_{C^1},
\|y_0-\mu^ny_n\Pi_n(s)\|_{C^1}\right\}\\
&\leq C\big(\mu^{3n}+n^3\mu^{4n}\big)
\leq C\mu^{3n}.
\end{align*}
The two differences on the first line define $\Delta_{1,n}$ and
$\Delta_{2,n}$. This finishes the proof of the lemma.
\end{proof}


\begin{lemma}
\label{lemma_In_estimate}
Recall that $I_n=x_0y_0$. There exists a constant $c>0$ such that
\begin{equation}
\label{eq_I_pure_epsilon}
I_n=c\mu^n-nac^2\mu^{2n}+O(\mu^{2n}).
\end{equation}
Moreover, $I_n-x_ny_n=O(\mu^{3n})$.
\end{lemma}

\begin{proof}
Let $s_n=\mu^nx_0y_n=O(\mu^n)$.
By Lemma~\ref{lemma_cross_form},
\begin{equation}
\label{eq_two_products}
I_n=x_0y_0=s_n\Pi_n(s_n)+O(\mu^{3n}),\qquad
x_ny_n=s_n\Pi_n(s_n)+O(\mu^{3n}),
\end{equation}
where recall
$$
\Pi_n(s)=1-nas+O(n^2s^2).
$$
In particular, $I_n=O(\mu^n)$ and
$x_ny_n-I_n=O(\mu^{3n})$.

Write $G=(G_1,G_2)$.
Since $p_0=G(p_n)$,
$$
G_2(x_n,y_n)=y_0=O(\mu^n).
$$
Homoclinic transversality gives
$\frac{\partial G_2}{\partial y}(p_{\rm out})\neq0$. Hence the implicit
function theorem gives $y_n=y_{\rm out}+O(\mu^n)$. Also
$x_0=G_1(x_n, y_n)=x_{\rm in}+O(\mu^n)$ because
$G(p_{\rm out})=p_{\rm in}=(x_{\rm in},0)$. Thus, if
$c=x_{\rm in}y_{\rm out}$, then
$$
x_0y_n=c+O(\mu^n)
$$
and therefore $s_n=c\mu^n+O(\mu^{2n})$. Substituting this into the first
equation in~\eqref{eq_two_products} gives
$$
I_n=c\mu^n-nac^2\mu^{2n}
+O(\mu^{2n}).
$$
\end{proof}

\begin{lemma}
\label{lemma_period_In}
There exists a
constant $d$ such that
\begin{equation}
\label{eq_period_In}
T_n=nT_0+T_{\rm exc}(p_{\rm out})+dI_n+nK(p,T_0)I_n
+O(n\mu^{2n}).
\end{equation}
\end{lemma}

\begin{proof}
Iterating forward in the first coordinate and backward in the second gives
(cf.~\cite[Lemma~4.1]{GRH})
$$
x_i=O(\mu^i),\qquad y_i=O(\mu^{n-i}),\qquad 0\le i\le n.
$$
Hence
$$
\sum_{i=0}^{n-1}(x_i+y_i)=O(1).
$$
Multiplying the two coordinates in~\eqref{eq_finite_moser} and using that
$\rho_1$ and $\rho_2$ are $C^1$ and vanish at the origin gives
\begin{equation}
\label{eq_product_drift}
x_{i+1}y_{i+1}-x_iy_i
=O\big((x_iy_i)^3(x_i+y_i)\big).
\end{equation}
Summing~\eqref{eq_product_drift} gives
\begin{equation}
\label{eq_product_almost_constant}
x_iy_i=I_n+O(\mu^{3n}),\qquad 0\le i\le n.
\end{equation}

Recall the splitting of $T_n$ into the local passage and the excursion part~\eqref{eq_split}.
Using~\eqref{eq_finite_roof} we can rewrite this decomposition as 
$$
T_n=\sum_{i=1}^n(T_0+\widehat \tau(x_i,y_i))+\widehat T_{\rm exc}(p_n),
$$
where $\widehat T_{\rm exc}(p_n)=T_{\rm exc}(p_n)-v(p_n)+v(p_0)$.
Again using using~\eqref{eq_finite_roof} and~\eqref{eq_product_almost_constant}, we obtain
\begin{equation}
\label{eq_sum_roof}
\sum_{i=0}^{n-1}\widehat\tau(x_i,y_i)
=nK(p,T_0)I_n+O(n\mu^{2n}).
\end{equation}
Hence
\begin{equation*}
\label{eq_asymp1}
T_n=nT_0+nK(p,T_0)I_n+O(n\mu^{2n})+\widehat T_{\rm exc}(p_n).
\end{equation*}
It remains to understand the remaining error $\widehat T_{\rm exc}(p_n)-\widehat T_{\rm exc}(p_{\rm out})=\widehat T_{\rm exc}(p_n)- T_{\rm exc}(p_{\rm out})$ (recall that $v$ vanishes on the coordinate axes). To calculate this error we first calculate, the local curve on which the successive points $p_n$ land (to the first order), and then apply the differential $D\widehat T_{\rm exc}(p_{\rm out})$ to appropriately scaled tangent vector to this curve. The same proof was used \cite[Lemma 4.3]{GRH} to calculate asymptotics of excursion time, where linearizing coordinates were used, which simplifies the calculation. Here we cannot use linearized coordinates because we need a better control on the error.

 To expand
$\widehat T_{\rm exc}(p_n)$ at $p_{\rm out}=(0,y_{\rm out})$, for small $\varepsilon>0$
define
\begin{equation*}
\label{eq_Phi_epsilon}
\Phi_\varepsilon(\xi,\eta)
=\big(\xi-G_1(\varepsilon\xi,\eta),
G_2(\varepsilon\xi,\eta)-\varepsilon\eta\big)
\end{equation*}
for $(\xi,\eta)$ near $(x_{\rm in},y_{\rm out})$. We have
$\Phi_0(x_{\rm in},y_{\rm out})=(0,0)$ and
$$
D\Phi_0(x_{\rm in},y_{\rm out})
=
\begin{pmatrix}
1&-\dfrac{\partial G_1}{\partial y}(p_{\rm out})\\[6pt]
0&\dfrac{\partial G_2}{\partial y}(p_{\rm out})
\end{pmatrix}.
$$
Its determinant is
$
\frac{\partial G_2}{\partial y}(p_{\rm out})
$
doesn't vanish by homoclinic transversality. Therefore, the implicit function theorem therefore gives
a local $C^{r-3}$ curve
$$
\varepsilon\longmapsto
(\xi(\varepsilon),\eta(\varepsilon)),\qquad
(\xi(0),\eta(0))=(x_{\rm in},y_{\rm out}),
$$
such that $\Phi_\varepsilon(\xi(\varepsilon),\eta(\varepsilon))=0$.

Now let $\varepsilon=\mu^n$. It follows from~\eqref{eq_cross_form}
and~\eqref{eq_passage_polynomial} that
$$
x_n=\varepsilon x_0+O(n\varepsilon^2),\qquad
y_0=\varepsilon y_n+O(n\varepsilon^2).
$$
Since $p_0=G(p_n)$, we have
\begin{align*}
\Phi_\varepsilon(x_0,y_n)
={}&\big(G_1(x_n,y_n)-G_1(\varepsilon x_0,y_n),\\
&G_2(\varepsilon x_0,y_n)-G_2(x_n,y_n)
+y_0-\varepsilon y_n\big)
=O(n\varepsilon^2).
\end{align*}
The local inverses in the implicit function theorem have uniformly bounded
derivatives. Hence
$$
(x_0,y_n)=(\xi(\varepsilon),\eta(\varepsilon))+O(n\varepsilon^2).
$$
Consequently,
$$
p_n=\big(\varepsilon\xi(\varepsilon),\eta(\varepsilon)\big)
+O(n\varepsilon^2).
$$
Differentiating with respect to $\varepsilon$ at zero gives
$$
p_n=p_{\rm out}
+\varepsilon\big(\xi(0),\eta'(0)\big)+O(n\varepsilon^2)
=p_{\rm out}
+\varepsilon\big(x_{\rm in},\eta'(0)\big)+O(n\varepsilon^2).
$$
By~\eqref{eq_I_pure_epsilon},
$I_n=c\varepsilon+O(n\varepsilon^2)$, where
$c=x_{\rm in}y_{\rm out}$. Therefore,
$$
p_n=p_{\rm out}
+\frac{I_n}{c}\big(x_{\rm in},\eta'(0)\big)+O(n\mu^{2n}).
$$
Let
$$
d=\frac1cD\widehat T_{\rm exc}(p_{\rm out})
\big(x_{\rm in},\eta'(0)\big).
$$
Taylor expansion of $\widehat T_{\rm exc}$ at $p_{\rm out}$ now gives
$$
\widehat{T}_{\rm exc}(p_n)
=\widehat{T}_{\rm exc}(p_{\rm out})+dI_n+O(n\mu^{2n})
=T_{\rm exc}(p_{\rm out})+dI_n+O(n\mu^{2n}),
$$
where, the last equality follows from
$v(p_{\rm in})=v(p_{\rm out})=0$.
Combining this with~\eqref{eq_sum_roof} yields
$$
T_n=nT_0+T_{\rm exc}(p_{\rm out})+dI_n+nK(p,T_0)I_n
+O(n\mu^{2n}).
$$
This proves the lemma.
\end{proof}

\begin{proof}[Proof of Proposition~\ref{prop_period_formula}]
Let $C=cd$. Substituting~\eqref{eq_I_pure_epsilon}
into~\eqref{eq_period_In} and using
$a=A(p,T_0)$ yields the posited formula
\begin{align*}
T_n={}&nT_0+T_{\rm exc}(p_{\rm out})+C\mu^n+cK(p,T_0)n\mu^n\\
&-c^2A(p,T_0)K(p,T_0)n^2\mu^{2n}
+O(n\mu^{2n}).
\end{align*}
\end{proof}

\section{Asymptotic formula for the multipliers}
\label{sec_multiplier_formula}

We continue to work with the shadowing setup of
Section~\ref{sec_period_formula}, using the refined Birkhoff coordinates.

\begin{proposition}
\label{prop_log_multiplier}
Let $\lambda_n>1$ be the unstable multiplier of $\gamma_n$. There exist
constants $C_0,C_1$ such that
\begin{equation}
\label{eq_log_multiplier}
\begin{split}
\log\lambda_n={}&-n\log\mu+C_0+C_1\mu^n+2cA(p,T_0)n\mu^n\\
&-\frac52c^2A(p,T_0)^2n^2\mu^{2n}+O(n\mu^{2n}).
\end{split}
\end{equation}
\end{proposition}

\begin{proof}
As before, let $s_n=\mu^nx_0y_n$. Let $\mathcal A_n$ be the derivative of
the polynomial map from Lemma~\ref{lemma_cross_form}:
\begin{multline}
\label{eq_cross_derivative}
\mathcal A_n=(\alpha_{ij})
=D_{(x_0,y_n)}
\big(\mu^nx_0\Pi_n(s_n),\mu^ny_n\Pi_n(s_n)\big)\\[6pt]
=\begin{pmatrix}
\mu^n\big(\Pi_n(s_n)+s_n\Pi_n'(s_n)\big)&
\mu^{2n}x_0^2\Pi_n'(s_n)\\[6pt]
\mu^{2n}y_n^2\Pi_n'(s_n)&
\mu^n\big(\Pi_n(s_n)+s_n\Pi_n'(s_n)\big)
\end{pmatrix}.
\end{multline}


The estimates~\eqref{eq_cross_F_F0_comparison} and
\eqref{eq_F0_polynomial_comparison} show that the derivative of
$(x_0,y_n)\mapsto(x_n,y_0)$ under $F^n$ is
$\mathcal A_n+O(\mu^{3n})$.
Moreover,
$$
\det\mathcal A_n
=\mu^{2n}\Pi_n(s_n)
\big(\Pi_n(s_n)+2s_n\Pi_n'(s_n)\big)
=O(\mu^{2n}).
$$
Since $\alpha_{22}=\mu^n(1+o(1))$, implicit differentiation applied to
$F^n(x_0,y_0)=(x_n,y_n)$ gives
\smallskip
\begin{equation}
\label{eq_scaled_local_derivative}
\alpha_{22}DF^n(x_0,y_0)
=\left(\begin{smallmatrix}
\det\mathcal A_n&\alpha_{12}\\-\alpha_{21}&1
\end{smallmatrix}\right)+O(\mu^{2n})
=\left(\begin{smallmatrix}
0&\alpha_{12}\\-\alpha_{21}&1
\end{smallmatrix}\right)+O(\mu^{2n}).
\end{equation}
\smallskip
Write $DG(p_n)=(\beta_{n,ij})$. Transversality of the homoclinic
intersection gives
$$
\frac{\partial G_2}{\partial y}(p_{\rm out})\neq0,\qquad\mathrm{and}
\qquad
\beta_{n,22}\longrightarrow
\frac{\partial G_2}{\partial y}(p_{\rm out}), n\to\infty.
$$
Multiplying~\eqref{eq_scaled_local_derivative} by $DG(p_n)$ gives
$$
\alpha_{22}DG(p_n)DF^n(p_0)=
\begin{pmatrix}
-\beta_{n,12}\alpha_{21}&
\beta_{n,11}\alpha_{12}+\beta_{n,12}\\
-\beta_{n,22}\alpha_{21}&
\beta_{n,21}\alpha_{12}+\beta_{n,22}
\end{pmatrix}
+O(\mu^{2n}).
$$
The main matrix on the right converges to
$$
\begin{pmatrix}
0&\dfrac{\partial G_1}{\partial y}(p_{\rm out})\\
0&\dfrac{\partial G_2}{\partial y}(p_{\rm out})
\end{pmatrix},
$$
whose non-zero eigenvalue we denote by
$$
P(0)=\frac{\partial G_2}{\partial y}(p_{\rm out}).
$$

Our goal is to derive the asymptotics as $n\to\infty$ of the unstable eigenvalue
$\lambda_n$. It is the larger eigenvalue of
$DG(p_n)DF^n(p_0)$. The eigenvalue of the scaled matrices
$\alpha_{22}DG(p_n)DF^n(p_0)$ converge: $\alpha_{22}\lambda_n\longrightarrow P(0)$, $n\to\infty$.
The corresponding eigenvector can be normalized as $(w_n,1)$, where $w_n$ is uniformly
bounded. 

The second coordinate of the eigenvalue equation gives
$$
\alpha_{22}\lambda_n
=\beta_{n,22}+\beta_{n,21}\alpha_{12}
-\beta_{n,22}\alpha_{21}w_n+O(\mu^{2n}).
$$
Note that, by~\eqref{eq_cross_derivative}, we have
$\alpha_{12}=O(n\mu^{2n})$ and $\alpha_{21}=O(n\mu^{2n})$. Since $w_n$ is uniformly bounded,
these terms can be absorbed into the error term, and hence
$$
\alpha_{22}\lambda_n=\beta_{n,22}+O(n\mu^{2n}).
$$
Multiplying the preceding identity by $\mu^n\Pi_n(s_n)/\alpha_{22}$ and
using~\eqref{eq_cross_derivative} gives
\begin{equation}
\label{eq_normalized_multiplier}
\mu^n\Pi_n(s_n)\lambda_n
=\frac{\Pi_n(s_n)}{\Pi_n(s_n)+s_n\Pi_n'(s_n)}
\beta_{n,22}+O(n\mu^{2n}).
\end{equation}
The main task is to expand $\beta_{n,22}$ and the coefficient in front of it. Recall from the proof of
Lemma~\ref{lemma_period_In} that $p_n$ lies approximately on the curve
$\varepsilon\mapsto(\varepsilon\xi(\varepsilon),\eta(\varepsilon))$ trough $p_{\rm out}$:
$$
p_n=\big(\mu^n\xi(\mu^n),\eta(\mu^n)\big)+O(n\mu^{2n}).
$$
Also recall that $G$ is $C^{r-3}$. Therefore,
$$
\beta_{n,22}
=\frac{\partial G_2}{\partial y}(p_n)
=\frac{\partial G_2}{\partial y}
\big(\mu^n\xi(\mu^n),\eta(\mu^n)\big)+O(n\mu^{2n}).
$$
Recall that by
Lemmas~\ref{lemma_cross_form} and~\ref{lemma_In_estimate},
\begin{equation}
\label{eq_In_sn}
I_n=s_n\Pi_n(s_n)+O(\mu^{3n}).
\end{equation}
Also, by~\eqref{eq_I_pure_epsilon},
$$
\mu^n=\frac{I_n}{c}+O(n\mu^{2n}).
$$
Define
$$
P(I)=\frac{\partial G_2}{\partial y}
\left(\frac{I}{c}\xi\left(\frac{I}{c}\right),
\eta\left(\frac{I}{c}\right)\right),
$$
which agrees with $P(0)$ value defined earlier.
The function $P$ is $C^2$, and the preceding expansion reads
$$
\beta_{n,22}=P(I_n)+O(n\mu^{2n}).
$$
Taylor expanding and using~\eqref{eq_passage_polynomial},
\eqref{eq_In_sn}, and~\eqref{eq_I_pure_epsilon}, gives

$$
\frac{\Pi_n(s_n)}{\Pi_n(s_n)+s_n\Pi_n'(s_n)}
=1+acn\mu^n-a^2c^2n^2\mu^{2n}+O(n\mu^{2n}).
$$

Plugging this expansion and the expansion for $\beta_{n,22}$ into~\eqref{eq_normalized_multiplier} we have
$$
\mu^n\Pi_n(s_n)\lambda_n
=P(I_n)\left(1+acn\mu^n-a^2c^2n^2\mu^{2n}\right)
+O(n\mu^{2n}).
$$
Taking logarithms and using~\eqref{eq_I_pure_epsilon} again gives
\begin{equation*}
\label{eq_log_normalized_multiplier}
\log\big(\Pi_n(s_n)\lambda_n\big)
=-n\log\mu+\log P(0)+c\frac{P'(0)}{P(0)}\mu^n
+acn\mu^n-\frac32a^2c^2n^2\mu^{2n}+O(n\mu^{2n}).
\end{equation*}

Recall that $s_n=c\mu^n+O(\mu^{2n})$.
Therefore, by~\eqref{eq_passage_polynomial},
$$
\log\Pi_n(s_n)
=-acn\mu^n+a^2c^2n^2\mu^{2n}+O(n\mu^{2n}).
$$
Subtracting this identity from the preceding expansion gives
$$
\log\lambda_n=-n\log\mu+\log P(0)+c\frac{P'(0)}{P(0)}\mu^n
+2acn\mu^n-\frac52a^2c^2n^2\mu^{2n}+O(n\mu^{2n}).
$$
Recalling that $a=A(p,T_0)$, this proves the posited formula~\eqref{eq_log_multiplier} with
$C_0=\log P(0)$, and $C_1=c {P'(0)}/{P(0)}$.
\end{proof}

\section{Proofs}

Given a periodic point $p=X^{T_0}(p)$ we have derived the asymptotic formulae for the periods on for the Lyapunov exponents:
$$
\begin{aligned}
T_n={}&nT_0+T_{\rm exc}(p_{\rm out})+C\mu^n\\
&+cK(p,T_0)n\mu^n
-c^2A(p,T_0)K(p,T_0)n^2\mu^{2n}+O(n\mu^{2n}),\\[3pt]
\log\lambda_n={}&-n\log\mu+C_0+C_1\mu^n\\
&+2cA(p,T_0)n\mu^n
-\frac52c^2A(p,T_0)^2n^2\mu^{2n}+O(n\mu^{2n}).
\end{aligned}
$$
If volume is the measure of maximal entropy then the synchronization
identity $\log\lambda_n=hT_n$ implies that, after multiplying the first equation by $h$, all coefficients in the above formulae must be equal. Coefficients by $n\mu^n$ give $hK(p, T_0)=2 A(p,T_0)$ while the coefficients by $n^2\mu^{2n}$ give $hK(p,T_0)=\frac52A(p,T_0)$. Hence $A(p,T_0)=K(p,T_0)=0$ for every periodic orbit.
Hence, by the Livshits theorem,
$[K]=0$ and Theorem~\ref{thm_Foulon_Hasselblatt} implies that $X^t$ is either a
contact flow or a constant-roof suspension.

In the contact case, Foulon's proof~\cite{Fou} applies and gives that $X^t$ is $C^{r_*}$ conjugate to a finite lift of
the constant-speed geodesic flow of a closed hyperbolic surface.\footnote{Foulon states his theorem for $C^\infty$ flows but the proof works for $C^2$ flows.}

In the suspension case, let $f$ be the return map to the torus section $\T^2$. Then $f$ is area preserving and area is the measure
of maximal entropy of $f$. Then de le Llave rigidity theorem implies that the continuous conjugacy to the linear model is, in fact, a $C^{r_*}$
regular~\cite{dlL}. Hence the proof of Theorem~\ref{thm_main} is complete.

\subsection{Proof of Corollary~\ref{cor_torus}}
\label{sec_torus_corollary}

Denote by $m$ the invariant area and let $\varphi=-\log J^u f$.
By the assumption, there is a $C^7$ function $\psi$ cohomologous to $\varphi$.
A sufficiently long finite Birkhoff
average of $\psi$ is negative. This average is still $C^7$ and cohomologous
to $\varphi$; call it $-\rho$. Then $\rho$ is a positive $C^7$ function.

The invariant area $m$ is the equilibrium state of the geometric potential
$\varphi$, and $P(\varphi)=0$. Since $-\rho$ is cohomologous to $\phi$,  $m$ is also the equilibrium state of
$-\rho$ and $P(-\rho)=0$. The description of equilibrium states for
suspensions~\cite[Proposition~6.1]{PP} implies that, for the
suspension flow over $f$ with roof $\rho$,
$$
\frac{m\times dt}{\int\rho\,dm}
$$
is its measure of maximal entropy. This is a smooth invariant measure, therefore,
Theorem~\ref{thm_main} applies.

A contact flow cannot admit a global invariant section. Thus it is a constant roof suspension and $f$ is smoothly conjugate to a linear hyperbolic automorphism.

\appendix

\section{Normal forms}
\label{app_normal_forms}

\subsection{Strong adapted charts}

\begin{proof}[Proof of Proposition~\ref{prop_nonstationary_birkhoff}]
Hurder and Katok~\cite[Proposition~4.2]{HK} constructed adapted charts with all needed properties
except properties~2 and~4. Property~2 is easy to add by pushing the charts along the flow lines to align the axes with strong stable and unstable manofolds. Thus it remains to establish property~4. 
Denote these initial adapted charts by $\widehat\Psi_q$ and their coordinate Poincar\'e
maps by $\widehat{\mathcal F}_{q,t}$.

Fix $T>0$ sufficiently large that
$$
\vartheta:=\sup_{q\in M}\mu(q,T)<1.
$$
For $q\in M$, write $\mu=\mu(q,T)$. Hurder and Katok
\cite[Lemma~4.7]{HK} give the following form:
$$
\widehat{\mathcal F}_{q,T}(x,y)
=\big(\mu x+\mu P_{q,T}(x,y),
\mu^{-1}y+\mu^{-1}Q_{q,T}(x,y)\big)
+R_{q,T}(x,y),
$$
where $P_{q,T}$ and $Q_{q,T}$ are homogeneous polynomials of degree $3$ and
$
D^jR_{q,T}(0,0)=0,\, 0\leq j\leq3.
$

Let us first determine the possible cubic terms. Axis invariance implies that
the polynomial $P_{q,T}$ is divisible by $x$ and $Q_{q,T}$ is divisible by $y$.
The terms of degree two in the identity
$\det D\widehat{\mathcal F}_{q,T}=1$ give
$$
\partial_xP_{q,T}+\partial_yQ_{q,T}=0.
$$
It follows that there are coefficients $\widehat\alpha=\widehat\alpha(q,T)$,
$A=A(q,T)$ and $\widehat\gamma=\widehat\gamma(q,T)$ such that
\begin{equation}
\label{eq_cubic_adapted_preliminary}
\big(P_{q,T}(x,y),Q_{q,T}(x,y)\big)=
\left(
\widehat\alpha x^3-Ax^2y+3\widehat\gamma xy^2,
-3\widehat\alpha x^2y+Axy^2-\widehat\gamma y^3
\right).
\end{equation}
We now eliminate $\widehat\alpha$ and $\widehat\gamma$ by adjusting the charts.

Given functions $\kappa,\nu$ on $M$, let $\Theta_q$ be the time-one map of the
Hamiltonian vector field generated by
$$
H_q(x,y)=\kappa(q)x^3y+\nu(q)xy^3.
$$
Explicitly, this vector field is
$$
X_{H_q}=(\partial_yH_q)\partial_x-(\partial_xH_q)\partial_y
=\big(\kappa(q)x^3+3\nu(q)xy^2\big)\partial_x
-\big(3\kappa(q)x^2y+\nu(q)y^3\big)\partial_y.
$$
The map $\Theta_q$ is an area-preserving local diffeomorphism tangent to the
identity and preserves both axes. Moreover,
\begin{equation}
\label{eq_cubic_chart_change}
\Theta_q(x,y)=(x,y)+U_q(x,y)+O\big((|x|+|y|)^5\big),
\end{equation}
where
$$
U_q(x,y)=
\left(\kappa(q)x^3+3\nu(q)xy^2,-3\kappa(q)x^2y-\nu(q)y^3\right).
$$
Set $\Psi_q=\widehat\Psi_q\circ\Theta_q$. The corresponding coordinate map is
$$
\mathcal F_{q,T}
=\Theta_{X^T(q)}^{-1}\circ\widehat{\mathcal F}_{q,T}\circ\Theta_q.
$$
Comparing cubic terms, the two non-resonant coefficients change according to
$$
\alpha(q,T)=\widehat\alpha(q,T)+\kappa(q)-\mu^2\kappa(X^T(q)),\qquad
\gamma(q,T)=\widehat\gamma(q,T)+\nu(q)-\mu^{-2}\nu(X^T(q)),
$$
while $A$ is unchanged.

To eliminate $\alpha(q,T)$, we solve
$$
0=\widehat\alpha(q,T)+\kappa(q)-\mu(q,T)^2\kappa(X^T(q)).
$$
Summing into the future gives
\begin{equation}
\label{eq_stable_cubic_adjustment}
\kappa(q)=-\widehat\alpha(q,T)
-\sum_{n=1}^{\infty}\mu(q,nT)^2\widehat\alpha(X^{nT}(q),T).
\end{equation}
The series converges uniformly because
$\mu(q,nT)\leq\vartheta^n$. After decreasing the H\"older exponent if necessary,
the convergence also holds in the H\"older norm, and hence its sum defines a
H\"older function $\kappa$ on $M$. The
$\kappa$-term removes the
$x^3$ term in the first component and, at the same time, the corresponding
$x^2y$ term in the second component.

Similarly, to eliminate $\gamma(q,T)$, we solve
$$
0=\widehat\gamma(q,T)+\nu(q)-\mu(q,T)^{-2}\nu(X^T(q)).
$$
Summing into the past gives
\begin{equation}
\label{eq_unstable_cubic_adjustment}
\nu(q)=\sum_{n=1}^{\infty}
\mu(X^{-nT}(q),nT)^2\widehat\gamma(X^{-nT}(q),T).
\end{equation}
This series also define a H\"older continuous function 
$\nu$ which removes the $y^3$ term in the second component together
with the corresponding $xy^2$ term in the first component, without
reintroducing the terms removed by the $\kappa$-term. 

It remains to pass from the fixed time $T$ to arbitrary time. Denote the two
non-resonant coefficients in the adjusted charts by $\alpha(q,t)$ and
$\gamma(q,t)$. Calculating the cubic terms in
$$
\mathcal F_{q,t+T}
=\mathcal F_{X^t(q),T}\circ\mathcal F_{q,t}
=\mathcal F_{X^T(q),t}\circ\mathcal F_{q,T}
$$
gives
$$
\alpha(q,t)=\mu(q,T)^2\alpha(X^T(q),t),
\qquad
\gamma(q,t)=\mu(q,T)^{-2}\gamma(X^T(q),t).
$$
Iterating the first identity into the future and the second one into the past
shows that
$$
\alpha(q,t)=\gamma(q,t)=0
$$
for every $q\in M$ and every $t>0$. Consequently,
$$
\mathcal F_{q,t}(x,y)
=\big(\mu(q,t)x(1-A(q,t)xy),
\mu(q,t)^{-1}y(1+A(q,t)xy)\big)+R_{q,t}(x,y).
$$
This is precisely~\eqref{eq_adapted_return}. 
\end{proof}

\subsection{The finite-order Birkhoff normal form}

The analytic normal form was proved by Moser~\cite{Mos56}. The first-order
finite-regularity normalization is carried out in
\cite{GS90}; see also~\cite{GS00}. Higher-order finite-regularity versions appeared
in~\cite{Gon13,DGG}. Neither reference gives precisely the statement
we need, because we require a $C^{r-3}$ regularity and remainders
$x^3y^2\rho_1$, $x^2y^3\rho_2$ with
$\rho_1,\rho_2\in C^{r-6}$ vanishing at the origin. We follow essentially the
same argument as in the above references, keeping track of the regularity and the remainders. After each
change of coordinates, we abuse notation and continue to denote the
conjugated return map by $F$.

\begin{proof}[Proof of Proposition~\ref{prop_normal_form}]$\,$\\
\noindent{\bfseries Step 1: Area preservation, invariant axes and their
linearization.}
Choose a $C^r$ transversal $\cT_p$ through $p$ which contains the local
strong stable and strong unstable manifolds of $p$. In local coordinates
$(u,v)$ on $\cT_p$, write
$$
(\iota_X\Omega)|_{\cT_p}=\rho(u,v)\,du\wedge dv.
$$
The function $\rho$ is $C^r$. After changing coordinates by
$$
\xi=u,\qquad \eta=\int_0^v\rho(u,s)\,ds,
$$
the transverse area is $d\xi\wedge d\eta$. We now write
the local stable manifold of $(0,0)$ as $\eta=s(\xi)$ and straighten it by the
area-preserving shear
$$
(\xi,\eta)\longmapsto(\xi,\eta-s(\xi)).
$$
We similarly straighten the local unstable manifold of $F$. These changes
do not lose any regularity. Hence the return map $F=(F_1,F_2)$ is $C^r$,
area-preserving and preserves both axes, and we will keep using $(x,y)$-variables.

We first record the elementary decomposition used at every step. If $h$ is
$C^2$, then
\begin{equation}
\label{eq_hadamard_decomposition}
h(x,y)=h(x,0)+h(0,y)-h(0,0)+
xy\int_0^1\int_0^1
\frac{\partial^2h}{\partial x\partial y}(tx,sy)\,dt\,ds.
\end{equation}
The one-dimensional linearization gives $C^r$
coordinates $\chi_s$ and $\chi_u$, such that
$$
\chi_s(F_1(x,0))=\mu\chi_s(x),\qquad
\chi_u(F_2(0,y))=\mu^{-1}\chi_u(y).
$$
The first coordinate is extended by
$$
(x,y)\longmapsto\left(\chi_s(x),\frac{y}{\chi_s'(x)}\right),
$$
and the second one by
$$
(x,y)\longmapsto\left(\frac{x}{\chi_u'(y)},\chi_u(y)\right).
$$
Both changes preserve area and both axes. Here we lose one derivative. Thus, in
the resulting $C^{r-1}$ coordinates,
$$
F_1(x,0)=\mu x,\qquad F_2(0,y)=\mu^{-1}y.
$$
Area preservation on the two axes also gives
$$
\partial_xF_1(0,y)=\mu,\qquad
\partial_yF_2(x,0)=\mu^{-1}.
$$
Therefore,
$$
F_1(x,y)=\mu x+x^2yh_1(x,y),\qquad
F_2(x,y)=\mu^{-1}y+xy^2h_2(x,y),
$$
where $h_1,h_2\in C^{r-4}$.  Applying
\eqref{eq_hadamard_decomposition} to $h_1$ and $h_2$ gives
\begin{align}
F_1(x,y)={}&\mu x\big(1+
(c_1+\varphi_1(x)+\psi_1(y))xy\big)
+x^3y^2R_1(x,y),
\label{eq_extended_order_one}\\
F_2(x,y)={}&\mu^{-1}y\big(1+
(c_2+\varphi_2(x)+\psi_2(y))xy\big)
+x^2y^3R_2(x,y).
\notag
\end{align}
where explicitely
$$
\begin{aligned}
c_1&=\mu^{-1}h_1(0,0),&
c_2&=\mu h_2(0,0),\\
\varphi_1(x)&=\mu^{-1}\big(h_1(x,0)-h_1(0,0)\big),&
\varphi_2(x)&=\mu\big(h_2(x,0)-h_2(0,0)\big),\\
\psi_1(y)&=\mu^{-1}\big(h_1(0,y)-h_1(0,0)\big),&
\psi_2(y)&=\mu\big(h_2(0,y)-h_2(0,0)\big).
\end{aligned}
$$
Hence $\varphi_j,\psi_j\in C^{r-4}$ and vanish at 0. Functions  $R_1$ and $R_2$ are $C^{r-6}$.

\medskip
\noindent{\bfseries Step 2: Removal of the non-resonant terms.}
Following~\cite[Section~8]{DGG}, we first remove $\varphi_1$. We seek a canonical change of coordinates $(x,y)\to (\xi,\eta)$ generated by
$$
V(x,\eta)=x\eta+\frac12(x\eta)^2\omega(x),
$$
where $\omega(0)=0$ is to be determined. Thus
$$
\xi=\partial_\eta V=x+x^2\eta\omega(x),\qquad
y=\partial_xV=\eta+\frac12\eta^2\frac{d}{dx}\big(x^2\omega(x)\big).
$$
Substituting this change in~\eqref{eq_extended_order_one} and comparing the
$x^2y$ terms replaces $\varphi_1(x)$ by
$$
\varphi_1(x)+\omega(\mu x)-\omega(x).
$$
Therefore the change removes $\varphi_1$ provided that
$$
\omega(\mu x)=\omega(x)-\varphi_1(x).
$$
This equation is solved directly by summing into the future
\begin{equation*}
\label{eq_weighted_solution}
\omega(x)=\sum_{k=0}^{\infty}\varphi_1(\mu^kx).
\end{equation*}
Indeed, the series converges in $C^{r-4}$ because $\varphi_1(0)=0$.
We need to have better regularity of $x^2\omega(x)$. By the definition of
$\varphi_1$,
$$
x^2\omega(x)=\sum_{k=0}^{\infty}\mu^{-2k}
\left(\mu^{-1}\partial_yF_1(\mu^kx,0)-c_1(\mu^kx)^2\right).
$$
It is easy to check that the series converges in $C^{r-2}$,
and $x^2\omega(x)$ is  $C^{r-2}$.
Solving $y=\partial_xV$ for $\eta$ gives the explicit coordinate change
\begin{equation}
\label{eq_canonical_change}
\eta(x,y)=\frac{2y}{1+\sqrt{1+2y\frac{d}{dx}(x^2\omega(x))}},\qquad
\xi(x,y)=x+\eta(x,y)x^2\omega(x).
\end{equation}
Since $x^2\omega(x)$ is $C^{r-2}$, this defines an area-preserving local $C^{r-3}$ diffeomorphism. By
construction, this change removes $\varphi_1$. The symmetric construction gives a canonical coordinate change $N_\psi$ which removes
$\psi_2$. Note that we do not construct the second canonical change after the first one. Instead, the second change is
constructed from the original $C^{r-1}$ axis-linearized map and hence is also $C^{r-3}$. In this way there is no additional loss of regularity. We explain why the two changes can be applied successively without altering the coefficient removed by the other. Let $N_\varphi$ denote the first change. Since $x^2\omega(x)$ has vanishing two-jet at the origin, we have 
$$
N_\varphi(0,y)=(0,y),\qquad DN_\varphi(0,y)=\operatorname{Id}.
$$
Let $F_\varphi=N_\varphi\circ F\circ N_\varphi^{-1}$. Since
$F(0,y)=(0,\mu^{-1}y)$,
$$
DF_\varphi(0,y)
=DN_\varphi(0,\mu^{-1}y)DF(0,y)DN_\varphi^{-1}(0,y)
=DF(0,y).
$$
Thus $\partial_x(F_\varphi)_2(0,y)=\partial_xF_2(0,y)$, so
$c_2+\psi_2(y)$ is unchanged. Similarly, $N_\psi$ and its derivative are the identity along the $x$-axis, so it does not reintroduce $\varphi_1$. Both
changes are $C^{r-3}$, and their composition is still $C^{r-3}$.

Now we examine the error terms after removing $\phi_1$ and $\psi_2$. Since both $\omega(x)/x$ and
$\omega'(x)$ are $C^{r-5}$, formula~\eqref{eq_canonical_change} gives
\begin{align*}
N_\varphi^{\pm1}(x,y)
&=\big(x+x^3y p_1^\pm(x,y),y+x^2y^2p_2^\pm(x,y)\big),\\
N_\psi^{\pm1}(x,y)
&=\big(x+x^2y^2q_1^\pm(x,y),y+xy^3q_2^\pm(x,y)\big),
\end{align*}
where $p_j^\pm,q_j^\pm\in C^{r-5}$ and $N_\psi$ denotes the symmetric
change. Direct substitution, together with the two homological equations,
gives
\begin{align*}
F_1(x,y)={}&\mu x\big(1+(c_1+\psi_1(y))xy\big)
+x^3y^2Q_1(x,y)+x^3yE_1(x,y),\\
F_2(x,y)={}&\mu^{-1}y\big(1+(c_2+\varphi_2(x))xy\big)
+x^2y^3Q_2(x,y)+xy^3E_2(x,y),
\end{align*}
where $Q_1,Q_2\in C^{r-6}$, $E_1,E_2\in C^{r-5}$, and
$
E_1(x,0)=0, E_2(0,y)=0.
$
Consequently,
$$
R_1(x,y)=Q_1(x,y)+\int_0^1\partial_yE_1(x,ty)\,dt,
\qquad
R_2(x,y)=Q_2(x,y)+\int_0^1\partial_xE_2(tx,y)\,dt
$$
are $C^{r-6}$ and
\begin{align*}
F_1(x,y)&=\mu x\big(1+(c_1+\psi_1(y))xy\big)+x^3y^2R_1(x,y),\\
F_2(x,y)&=\mu^{-1}y\big(1+(c_2+\varphi_2(x))xy\big)
+x^2y^3R_2(x,y).
\end{align*}
Area preservation, $\det DF=1$, immediately gives
$$
c_2=-c_1,\qquad \varphi_2=\psi_1=0.
$$


\medskip
\noindent{\bfseries Step 3: The second Birkhoff coefficient.}
Let $b=\mu^{-1}R_1(0,0).$
The coefficient of $x^2y^2$ in $\det DF=1$ gives
$$
3\big(b+\mu R_2(0,0)+c_1c_2\big)=0.
$$
Thus $\mu R_2(0,0)=c_1^2-b$. It remains to set $a=-c_1$ and
$$
\rho_1(x,y):=R_1(x,y)-R_1(0,0),\qquad
\rho_2(x,y):=R_2(x,y)-R_2(0,0).
$$
to obtain the posited Birkhoff normal form~\eqref{eq_finite_moser}.

\medskip
\noindent{\bfseries Step 4: The return time.}
Let $\Psi$ be the resulting $C^{r-3}$ parametrization of $\cT_p$, and let
$T_0+\tau$ be the return time. Since $\cT_p$ contains the local strong stable
and strong unstable manifolds we have
$$
\tau(x,0)=\tau(0,y)=0.
$$
Moreover, $\tau\in C^{r-3}$. Dividing successively by
$x$ and $y$ we get
$$
\tau(x,y)=xy\,q(x,y),
$$
where $q$ is $C^{r-5}$ and hence at least $C^2$.
By chart independence of the periodic value of the longitudinal cocycle
\cite[Lemma~2.2]{GRH},
$$
K(p,T_0)=\partial_{xy}\tau(0,0)=q(0,0).
$$

Applying~\eqref{eq_hadamard_decomposition} to $q$ gives
$$
q(x,y)=q(x,0)+q(0,y)-K(p,T_0)+O(xy).
$$
After multiplication by $xy$, the two nonconstant axis terms are precisely
the terms beginning in degree three which are not already $O((xy)^2)$. We
will remove them by adding a coboundary $v\circ F-v$, where $v$ has the
form $v=xyw$. Define
$$
w(x,y)=\sum_{j=0}^{\infty}\big(q(\mu^jx,0)-K(p,T_0)\big)
+\sum_{j=1}^{\infty}\big(K(p,T_0)-q(0,\mu^jy)\big).
$$
The series converge in $C^{r-5}$, and
$$
w(\mu x,\mu^{-1}y)-w(x,y)
=2K(p,T_0)-q(x,0)-q(0,y).
$$
Set $\widehat\tau=\tau+v\circ F-v$. Then
\begin{align*}
\big|\widehat\tau(x,y)-K(p,T_0)xy\big|
\le{}&|xy|\big|q(x,y)+w(\mu x,\mu^{-1}y)-w(x,y)-K(p,T_0)\big|\\
&+\big|F_1(x,y)F_2(x,y)w(F(x,y))
-xyw(\mu x,\mu^{-1}y)\big|.
\end{align*}
The first term on the right is $O((xy)^2)$ by the preceding discussion.
It remains to estimate the second term. By the triangle inequality,
\begin{multline*}
\big|F_1(x,y)F_2(x,y)w(F(x,y))
-xyw(\mu x,\mu^{-1}y)\big|\\
\le \big|F_1(x,y)F_2(x,y)-xy\big|
\big|w(F(x,y))\big|\\
+|xy|\big|w(F(x,y))-w(\mu x,\mu^{-1}y)\big|.
\end{multline*}
The normal form
\eqref{eq_finite_moser} gives $F_1(x,y)F_2(x,y)=xy+O((xy)^3)$, which takes care of the first term.
Since $w$ is $C^1$, we have
$$
\big|w(F(x,y))-w(\mu x,\mu^{-1}y)\big|
=O(x^2|y|+|x|y^2)=O(|xy|),
$$
Hence the second term in the first
inequality is $O((xy)^2)$, and
$$
\widehat\tau(x,y)=K(p,T_0)xy+O((xy)^2).
$$
\end{proof}

\end{document}